\documentclass[a4paper,11pt,onecolumn]{amsart}
\usepackage{amsfonts}
\usepackage{pb-diagram}
\title[Finitely generated nilpotent groups]{On some of the residual properties of finitely generated nilpotent groups}

\newtheorem{thm}{Theorem}
\newtheorem{lemma}[thm]{Lemma}

\newtheorem{cor}[thm]{Corollary}

\numberwithin{thm}{section}

\newcommand{\gam}{\Gamma}

\newcommand{\bZ}{\mathbb{Z}}

\newcommand{\bC}{\mathbb{C}}

\DeclareMathOperator{\rk}{rk}

\newcommand{\bQ}{\mathbb{Q}}

\author[T. Koberda]{Thomas Koberda}
\address{Department of Mathematics\\ Harvard University\\ 1 Oxford St.\\ Cambridge, MA 02138 }
\email{ koberda@math.harvard.edu}
\begin{document}
\begin{abstract}
In recent years, the RFRS condition has been used to analyze virtual fibering in $3$-manifold topology.  Agol's work shows that any $3$-manifold with zero Euler characteristic satisfying the RFRS condition on its fundamental group virtually fibers over the circle.  In this note we will show that a finitely generated nilpotent group is either virtually abelian or is not virtually RFRS, a result which may be of independent interest though not directly applicable to $3$-manifold topology.  As a corollary, we deduce that any RFRS group cannot contain a nonabelian torsion-free nilpotent group.  This result also illustrates some of the interplay between residual torsion-free nilpotence and the RFRS condition.
\end{abstract}
\maketitle
\begin{center}
\today
\end{center}

\section{Introduction}
Let $\gam$ be a finite undirected graph.  By this we mean that $\gam$ consists of a finite set $V$ of vertices and a finite set $E\subset (V\times V\setminus\Delta)/\bZ/2\bZ$ of edges.  We consider edges as unordered pairs of vertices so that the edges are undirected, and we assume that there is at most one edge connecting any two vertices.  Removing the diagonal from $V\times V$ prevents any edges looping at a vertex.  We define the {\bf graph group} $G(\gam)$ associated to $\gam$ by the presentation \[\langle V\mid [v_1,v_2] \textrm{ if $\{v_1,v_2\}\in E$}\rangle.\]  The graph group $G(\gam)$ is often called a right-angled Artin group (RAAG).

Recall that a group $G$ is called {\bf residually finite rationally solvable} or {\bf RFRS} if there exists an exhaustive filtration of normal (in $G$) finite index subgroups of $G$ \[G=G_0>G_1>G_2>\cdots\] such that $G_{i+1}>ker\{G_i\to G_i^{ab}\otimes\bQ\}$.  Agol proves:
\begin{thm}[\cite{A}]
Let $M$ be a $3$-manifold with $\chi(M)=0$, and suppose that $\pi_1(M)$ is RFRS.  Let $\phi\in H^1(M,\bZ)$ be a non-fibered cohomology class.  Then there is a finite cover $M'\to M$ such that the pullback of $\phi$ to $H^1(M',\bZ)$ of $\phi$ lies in the cone over the boundary of a fibered face of $H^1(M',\bZ)$.
\end{thm}

He then deduces:
\begin{cor}
Let $M$ be a $3$-manifold with $\chi(M)=0$ and suppose that $\pi_1(M)$ virtually embeds in a graph group.  Then $M$ virtually fibers over the circle.
\end{cor}

The interest in graph groups as they relate to virtual fibering is partially due to the following theorem of Wise (cf. the work in \cite{HW}):
\begin{thm}[\cite{W1},\cite{W2}]
Let $G$ be a word-hyperbolic group with a quasiconvex hierarchy.  Then $G$ has a finite index subgroup $G'$ which embeds in a graph group $R$.
\end{thm}
Wise uses Agol's theorem to deduce:
\begin{cor}
A Haken hyperbolic $3$-manifold of finite volume is virtually fibered.
\end{cor}

In \cite{A}, Agol provides a proof that right-angled reflection groups are virtually RFRS.  On the other hand:
\begin{thm}[Droms (\cite{D}), Duchamp and Thibon (\cite{DT}), A. Vijayan (\cite{C}), and others]\label{t:rtfn}
Let $G$ be a graph group.  Then $G$ is residually torsion-free nilpotent.
\end{thm}

Residual torsion-free nilpotence and the RFRS condition are independent since it is known that there exist residually torsion-free nilpotent groups which are not even virtually RFRS (see \cite{K}).  Examples of such groups are compact $3$-manifolds with Nil geometry and nontrivial circle bundles with trivial monodromy over a closed surface.

To further illustrate the incompatibility of RFRS and residual torsion-free nilpotence, we have the following result:

\begin{thm}\label{p:rtfn}
Let $N$ be a finitely generated nilpotent group.  Then $N$ is either virtually abelian or $N$ is not virtually RFRS.
\end{thm}

The RFRS condition was originally developed to answer questions about $3$-manifold groups, and nilpotent groups never arise as fundamental groups of finite volume hyperbolic $3$-manifolds.  Theorem \ref{p:rtfn} is not so much a result with applicability to $3$-manifold theory but rather one which might be of independent interest.

Since the RFRS condition is inherited by subgroups, we deduce:
\begin{cor}
Let $G$ be a virtually RFRS group.  Then $G$ does not contain a torsion-free nilpotent group.
\end{cor}

In particular it follows that graph groups do not contain torsion-free nilpotent groups, though this fact can be deduced from the Tits alternative for graph groups.

\section{Acknowledgments}
The author thanks Ruth Charney, Stefan Friedl and Dani Wise for useful conversations.

\section{Proof of the main result}
In \cite{K}, the author proves that the Heisenberg group is not RFRS.  An easy modification of the proof shows that the Heisenberg group is in fact not virtually RFRS.  The author also proves that a nontrivial circle bundle with trivial monodromy over a closed orientable surface has a residually torsion-free nilpotent fundamental group, but remarks that it cannot be RFRS by Agol's theorem.  The author then proves directly that the fundamental group of such a manifold cannot be RFRS.  We now adapt this discussion to general finitely generated nilpotent groups.

\begin{lemma}\label{l:vab}
Let $N$ be a finitely generated torsion-free nilpotent group which is virtually abelian.  Then $N$ is abelian.
\end{lemma}
\begin{proof}
Let $N'<N$ be a finite index normal subgroup of $N$ which is abelian.  Since $N'$ is finitely generated, it is isomorphic to $\bZ^n$, and the quotient $N/N'$ is some finite group $F$.  The conjugation action of $N$ on $N'$ factors through $F$, since $N'$ is abelian and therefore acts trivially.  We thus obtain a representation $F\to GL_n(\bZ)$, whose image is nontrivial since we may assume $N$ is nonabelian.  But the image of any element of $F$ in $GL_n(\bZ)$ is semisimple since $F$ is finite (it can be diagonalized over $\bC)$.  It follows that the action of $F$ on $N'$ cannot be unipotent.  In particular, let $f\in F$ have image $A_f\in GL_n(\bZ)$.  Then $A_f-I$ is not nilpotent, so that $(A_f-I)^M(\bZ^n)\neq \{0\}$ for all $M>0$.  It follows that $N$ could not have been nilpotent.
\end{proof}

What the proof of Lemma \ref{l:vab} shows is that if $A$ is a torsion-free abelian normal subgroup of a nilpotent group, then the conjugation action of $N$ on $A$ factors through $N/A$ and must act as a unipotent group of automorphisms of $A$, which is not possible.  If $T$ is a non-unipotent matrix, then putting $T$ in Jordan canonical form (over $\bC$) we see that $(T-I)^n$ has positive rank (even over $\bQ$), no matter how large $n$ gets.

Recall that finitely generated nilpotent groups have a well-defined notion of {\bf rank} (cf. \cite{R}).  Precisely, let $N$ be finitely generated and nilpotent.  Suppose that \[N=N_1^1>N_2^1>\cdots>N_i^1>N_{i+1}^1=\{1\}\] and \[N=N_1^2>N_2^2>\cdots>N_j^2>N_{j+1}^2=\{1\}\] are two normal filtrations of $N$ such that $N_i^1$ and $N_j^2$ are abelian and each successive quotient is abelian.  Then the rank of $N$, written $\rk N$, is given by \[\sum_{k=1}^i \rk N_k^1/N_{k+1}^1=\sum_{k=1}^j \rk N_k^2/N_{k+1}^2.\]  Such filtrations always exist, so that the rank is a well-defined nonnegative integer.  Notice that by definition, if \[1\to N'\to N\to N''\to 1\] is a short exact sequence of nilpotent groups, where $N$ is finitely generated, then $\rk N= \rk N'+\rk N''$.

Let $N$ be a torsion-free nonabelian nilpotent group with center $Z$.  It is standard that $Z$ is again finitely generated (see \cite{R}).  Suppose that the extension \[1\to Z\to N\to N/Z\to 1\] does not split.  Let $N'$ be a finite index subgroup of $N$.  $N'$ is still clearly torsion-free and nilpotent.  It is conceivable that the extension \[1\to Z\cap N'\to N'\to N'/(Z\cap N')\to 1\] does split.  However, as a consequence of Lemma \ref{l:vab} and by induction on the rank of $N$, we see that the extension \[1\to Z(N')\to N'\to N'/Z(N')\to 1\] does not split.

\begin{lemma}\label{l:virt}
Let $N$ be a finitely generated torsion-free nilpotent group and let $\{1\}\neq Z<Z(N)$ be a central subgroup.  Suppose that the extension \[1\to Z\to N\to N/Z\to 1\] is not a split extension.  Then either the restriction of $N\to N^{ab}$ to $Z$ is not injective or $N$ virtually splits as a nontrivial direct product.
\end{lemma}
\begin{proof}
By \cite{E}, we may suppose that $N/Z$ is torsion-free.  Note that it is not possible for $N$ to be abelian, for if it was then the extension would automatically split.  Suppose that the abelianization map restricted to $Z$ is injective.  Then we see that $Z$ virtually splits off as a direct summand of $N^{ab}$, so that $N^{ab}\cong A\times Z'$ for some abelian group $A$ and a torsion-free abelian group $Z'$ which contains the image of $Z$ with finite index.  Let $K$ be the kernel of the map $N\to Z'$.  Then $K$ fits into a (not necessarily central) extension of the form \[1\to K\to N\to A\times Z'\to 1.\]  This extension splits partially, since $Z$ is a subgroup of $N$.  Since $Z$ is central, $Z$ and $K$ normalize each other in $N$.  Since $Z$ maps injectively to $Z'$, $Z\cap K=\{1\}$.  It follows that $Z$ and $K$ together generate a subgroup of $N$ which is isomorphic to $K\times Z$.  Since $N/(K\times Z)\cong Z'/Z$, the conclusion follows.
\end{proof}

\begin{cor}\label{c:notinj}
Let $N$ be a finitely generated torsion-free nilpotent group, and suppose that $Z=Z(N)$ maps injectively to $N^{ab}$.  Then $N$ is abelian.
\end{cor}
\begin{proof}
By Lemma \ref{l:virt}, $N$ virtually splits as a direct product of the form $K\times Z$.  Since $K$ is nilpotent, it has a center which is clearly a subgroup of $N$ and is hence torsion-free, abelian, and finitely generated.  Since $Z$ is central, $K\times Z$ acts trivially on the center of $K$.  Since $Z$ is precisely the center of $N$, $X=Z(K)$ is not central in $N$.  It follows that $\gam=N/(K\times Z)$ acts nontrivially on $X$.  As $X\cong \bZ^n$, the conjugation action of $\gam$ is given by a representation $\gam\to GL_n(\bZ)$.  Since $\gam$ is a finite group (it is abelian as an easy consequence of the proof of Lemma \ref{l:virt}), $\gam$ acts by nontrivial semisimple automorphisms (they can in fact be diagonalized over $\bC$).  It follows that the conjugation action of $\gam$ is not unipotent, which violates the nilpotence of $N$.  It follows that if $K$ is nontrivial then $Z$ is not the whole center of $N$.  It follows that $N$ is virtually abelian and hence abelian.
\end{proof}

\begin{proof}[Proof of Theorem \ref{p:rtfn}]
Since $N$ is polycyclic, we see that $N$ is virtually torsion-free by \cite{R}.  So, we may assume that $N$ is a nonabelian torsion-free nilpotent group.  By the corollary to Lemma \ref{l:vab}, it suffices to show that $N$ is not RFRS.  Let $Z=Z(N)$ and let $\{G_i\}$ be any filtration of $N$ which witnesses the claim that $N$ is RFRS, with $G_0=N$.  By Corollary \ref{c:notinj}, $Z$ does not map injectively into $N^{ab}$.  There is an element $z\in Z$ such that $H=\langle z\rangle$ is infinite cyclic  and the canonical map $H\to N^{ab}$ is not injective.  The fact that $H\to N^{ab}$ is not injective implies that $H$ is contained in $ker\{N\to N^{ab}\otimes\bQ\}$, so that $H$ is contained in $G_1$.

We may assume that $H<G_i$ for some $i\geq 1$.  Since $G_i$ has finite index in $N$, we may assume that $G_i$ is not abelian.  Suppose that $H$ maps injectively to $G_i^{ab}$.  Then by Lemma \ref{l:virt}, $G_i$ virtually splits as a nontrivial direct product.  By Corollary \ref{c:notinj}, there is another central cyclic subgroup $H'$ of $G_i$ which does not map injectively to $G_i^{ab}$ and hence must be contained in $G_{i+1}$.

If \[\bigcap_i G_i=\{1\},\] then eventually each element $z$ of the center of $N$ must be contained in $G_i\setminus G_{i+1}$ for some index $i$.  If $i=i(z)$ is such an index then $G_i$ virtually splits as $K\times\bZ$, where the $\bZ$ factor generated by some power of $z$.  Note that $\rk K+1\leq \rk N$.  Since $K$ is nilpotent, it also has a center.  By induction on rank, $N$ is virtually abelian.  The conclusion follows by Lemma \ref{l:vab}.
\end{proof}

For completeness, we include the following observation, which is implicit in the work of Agol and others:
\begin{lemma}
Let $G$ be a group which is RFRS and $H<G$ a subgroup.  Then $H$ is RFRS.
\end{lemma}
\begin{proof}
Let $\{G_i\}$ be a filtration of $G$ which witnesses the fact that $G$ is RFRS.  Clearly $\{H_i\}=\{G_i\cap H\}$ is a filtration on $H$.  Since there is a map $H_i^{ab}\to G_i^{ab}$ given by the inclusion, if $h\in ker\{H_i\to H_i^{ab}\otimes\bQ\}$ then $h$ maps to a torsion element in $G_i^{ab}$ and is therefore contained in $G_{i+1}$.  It follows that $h$ is contained in $H_{i+1}$.  Conversely, if $h\in G_i\setminus G_{i+1}$ then $h$ maps to a non-torsion element of $G_i^{ab}$ and must therefore map to a non-torsion element of $H_i^{ab}$.
\end{proof}

It follows that if $G$ is a virtually RFRS group then $G$ cannot contain a nonabelian torsion-free nilpotent group.

\end{document}